\documentstyle[amssymb,euscript]{amsart}
\title{Filters and games.}
\author{Tomek Bartoszynski and Marion Scheepers}
\address{Department of Mathematics\\ Boise State University\\ Boise\\
   Idaho 83725}
\email{\tt tomek@@math.idbsu.edu, marion@@math.idbsu.edu}
\thanks{The second author was supported by Idaho State Board of
  Education grant 94--051.}
\newtheorem{theorem}{Theorem}

\newtheorem{definition}{Definition}

\newcommand{\naturals}{\Bbb N}

\newcommand{\G}{{\frak G}_2({\cal F})}
\newcommand{\HH}{{\frak G}_1({\cal F})}
\begin{document}

\maketitle
\begin{abstract}  We obtain game--theoretic characterizations for
  meagerness and rareness of filters on $\omega$.
\end{abstract}

   One of the classical methods for obtaining a set of real numbers
   which does not have the property
   of Baire, is to interpret appropriate filters on
   $\naturals=\{1,2,3,\dots\}$ as subsets of $[0,1]$. Filters which
   result in a set having the property of Baire have
   nice combinatorial characterizations, due to Talagrand \cite{T}:
   \begin{theorem}[Talagrand]
     For a filter ${\cal F}$ on $\omega$ the following are equivalent:
     \begin{enumerate}
       \item{${\cal F}$ does not have the property of Baire,}
       \item{The set of enumeration functions of sets in ${\cal F}$ is
           unbounded in $^{\omega}\omega$, ordered by eventual
           domination,}
       \item{For every partition of $\omega$ into disjoint finite
           sets, there is a set in ${\cal F}$ which is disjoint from
           infinitely many blocks in the partition.}
     \end{enumerate}
   \end{theorem}
   In particular, we see that a filter is meager if, and only if, it
   has the property of Baire. We use this without further notice below.
   
   These combinatorial characterizations suggest certain infinite
   two--player games. We study such a game in section 1. We prove that
   having the property of Baire is equivalent with the
   assertion that player ONE of our game has a winning strategy
   (Theorem \ref{meagerandONE}). In the second section we study a
   natural variation of the  first game, and show that here, being
   non-rare is equivalent to ONE having a winning strategy in
   this game (Theorem \ref{rareandONE}). 

   Our notation is mostly standard; the only exception may be that we
   use $\frown$ to denote concatenation of sequences.   
   From now on, let ${\cal
   F}\subset{\EuScript P}(\naturals)$ be a non--principal filter.

\section{The game $\HH$.}

   This game is played as follows: In the $k$-th inning, player ONE
   chooses $m_k\in\naturals$ and player TWO responds with
   $n_k\in\naturals$. TWO wins the play $(m_1,n_1,\dots,m_k,n_k,\dots)$
   if:
\begin{enumerate}
\item{$n_1<n_2<\dots<n_k<\dots$ and}
\item{there are infinitely many $k$ such that $m_k<n_k$ and}
\item{$\{n_1,n_2,\dots,n_k,\dots\}\in{\cal F}$.}
\end{enumerate}
   Otherwise, ONE wins.

\begin{theorem}\label{th1}
   TWO does not have a winning strategy in $\HH$.
\end{theorem}

\begin{pf} Suppose on the contrary that $F$ is a winning strategy for
TWO in $\HH$. Here, and below, we may assume without loss of
generality that $F$ depends on only ONE's moves -- TWO's moves cane be
decoded from these, using $F$.\\
{\bf Claim:} For each finite sequence $\sigma$ of natural numbers,
there exists a finite sequence $\tau$ of natural numbers and a natural
number $n$ such that
$F(\sigma\frown\tau\frown(y))>y$ for each $y>n$.\\
{\bf Proof of the claim:} If the claim were false, then there would be a
finite sequence $\sigma_0$ such that for every finite sequence $\tau$ 
and every natural number $n$, there is a $y>n$ such that
$F(\sigma_0\frown\tau\frown(y))\leq y$.
Fix $\sigma_0$ and choose $y_1<y_2<y_n<\dots$ so that
$max(\sigma_0)<y_1$ and
$F(\sigma_0\frown(y_1,\dots,y_n)\frown(y_{n+1}))\leq y_{n+1}$ for each
$n$. Then TWO lost the play where ONE begins the game by picking the
values of $\sigma_0$ first, and then the numbers $y_n$ consecutively.
This contradicts our hypothesis that $F$ is a winning strategy for
TWO.\\

Using the claim, we now derive a contradiction from the hypothesis
that $F$ is a winning strategy for TWO. To begin, let $\sigma$ be the
empty sequence. Pick $\tau_{\emptyset}$ as in the claim. Then choose
$y_0>\max\tau_{\emptyset}$ also as in the claim. Put
$\sigma_0=\tau_{\emptyset}\frown(y_0)$ and pick $\tau_0$ as in the
claim.

Next choose $y_1>\max\{F(\sigma_0\frown\tau_0\restriction_{j}):j<\omega\}$ and put
$\sigma_1=\tau_{\emptyset}\frown(y_1)$; choose $\tau_1$ as in the
claim.

Next choose $y_2>\max\{F(\sigma_1\frown\tau_1\restriction_{j}):j<\omega\}$ and put
$\sigma_2=\sigma_0\frown\tau_0\frown(y_2)$ and then choose $\tau_2$ as
in the claim.

Then choose $y_3>\max\{F(\sigma_2\frown\tau_2\restriction_{j}):j<\omega\}$
and put $\sigma_3=\sigma_1\frown\tau_1\frown(y_3)$; then choose
$\tau_3$ as in the claim.

Continuing like this we choose three sequences
\[(y_0,y_1,y_2,\dots),\]
\[(\sigma_0,\sigma_1,\sigma_2,\dots) \mbox{ and}\]
\[(\tau_0,\tau_1,\tau_2,\dots)\]
such that
\begin{enumerate}
\item{$\sigma_{2i+1}=\sigma_{2i-1}\frown\tau_{2
   i-1}\frown(y_{2i +1})$,}
\item{$\sigma_{2i+2}=\sigma_{2i}\frown\tau_{2i}\frown(y_{2i+2})$,}
\item{$\tau_{2i+1}$ is chosen as in the claim for
$\sigma_{2i+1}$,}
\item{$\tau_{2i}$ is chosen as in the claim for $\sigma_{2i}$,}
\item{$y_{2i+1}$ is chosen as in the claim for $\sigma_{2
   i-1}\frown \tau_{2i-1}$, so that $y_{2
   i+1}>\max\{F(\sigma_{2i}\restriction_j):j<\omega\}$ and}
\item{$y_{2i+2}$ is chosen as in the claim for $\sigma_{2
   i}\frown \tau_{2i}$, so that
   $y_{2i+2}>\max\{F(\sigma_{2i+1}\restriction_j):j<\omega\}$.} 
\end{enumerate}

   Let $f$ and $g$ be the unique sequences of natural numbers such that
   $\sigma_{2i}\subset f$ and $\sigma_{2i +1}\subset g$ for
   all $i$. Then by construction the sets 
\[A_f=\{F(f\restriction_n):n\in\naturals\}\mbox{ and}\]
\[A_g=\{F(g\restriction_n):n\in\naturals\}\]
have finite intersection, and both are response sets for TWO using the
strategy $F$. But then at least one of these sets is not in the
non--principal filter ${\cal F}$! Thus there is a play for ONE against
$F$ which defeats TWO.
\end{pf}

\begin{theorem}\label{meagerandONE} The following statements are equivalent:
\begin{enumerate}
\item{${\cal F}$ is a meager filter.}
\item{ONE has a winning strategy in $\HH$.}
\item{$\HH$ is determined.}
\end{enumerate}
\end {theorem}

\begin{pf} The implication $1\Rightarrow 2$ follows from the negation
  of $(2)$ in Talagrand's theorem. The implication $2\Rightarrow 3$ is trivial.
  We show that $(3)\Rightarrow(1)$. 

   We already know that TWO does not have a winning strategy.
   Let $F$ be a winning strategy
   for ONE in $\HH$. We may assume that $\max\{x_1,\dots,x_n\}+1\leq
   F(x_1,\dots,x_n)$ for all $x_1<\dots<x_n$. 
   Define a function $f$ so that $f(0)=F(\emptyset)$ and for $n>1$,
\[f(n)=\max\{F(t_1,t_2,\dots,t_j): t_i\leq n \mbox{ and }t_1<\dots<t_j
   \mbox{ for }i\leq j\leq n\}.\]
   Then $f$ is monotonic -- i.e., $f(m)\leq f(n)$ whenever $m\leq n$.
   For each $n$ define $h_n(m)=f^m(n)$ for each $m$. Then choose $g$ 
   so that $g(k)=\max\{h_m(n):m,n\leq k\}+1$ for all $k$. 

   We claim that if $X\in{\cal F}$, then $g$ dominates $enum_X$, the
   enumeration function of $X$. Fix such an $X$ and write
   $X=\{x_1,x_2,\dots\}$ in increasing order. Since $F$ is a winning
   strategy for ONE, the play
\[F(\emptyset), x_1, F(x_1), x_2, F(x_1,x_2), x_3,
F(x_1,x_2,x_3),\dots\]
   is won by ONE. Thus, there is an $N\in\naturals$ such that
   $x_{n+1}\leq F(x_1,\dots,x_n)$ for $n\geq N$.

   We see in particular that:\\
   $x_{N+1}\leq F(x_1,\dots,x_{N})\leq f(x_N)$,\\
   $x_{N+2}\leq F(x_1,\dots,x_{N+1})\leq f(x_{N+1})\leq f^2(x_N)$,\\
   and in general, $x_{N+k}\leq F(x_1,\dots, x_{N+k-1})\leq
   f(x_{N+k-1})\leq f^k(x_N) = h_{x_N}(k)$ for each $k$.

   But $g$ eventually dominates $h_{x_N}$. Choose $K>x_N$ so large
   that $h_{x_N}(k)<g(k)$ for all $k\geq K$. Then $x_{N+k}< g(k)$ for
   all $k\geq K$; in particular, $x_k< g(k)$ for all $k\geq K$.
   
   We have shown that the set of enumeration functions of elements of
   ${\cal F}$ is bounded; by Talagrand's theorem, ${\cal F}$ is meager.
\end{pf}

\section{The game $\G$.}

   In the game $\HH$, TWO's objective was to play a sequence
   (enumerating a set from the filter) not
   eventually dominated by ONE's sequence. What is the situation when
   we change TWO's objective to playing a sequence (enumerating  a set
   from the filter) which actually eventually dominates ONE's
   sequence? We consider this now: In the $k$-th inning, player ONE
   chooses $m_k\in\naturals$ and player TWO responds with
   $n_k\in\naturals$. TWO wins the play $(m_1,n_1,\dots,m_k,n_k,\dots)$
   if:
\begin{enumerate}
\item{$(n_1,\dots,n_k,\dots)$ eventually dominates
   $(m_1,m_2,\dots,m_k,\dots)$ \underline{and}}
\item{$\{n_1,n_2,\dots,n_k,\dots\}\in{\cal F}$.}
\end{enumerate}
   Otherwise, ONE wins.

\begin{theorem} TWO does not have a winning strategy in $\G$.
\end{theorem}

\begin{pf} This theorem follows immediately from Theorem \ref{th1},
  because the game $\G$ is harder for TWO than $\HH$. However, there
  is a much simpler argument than that for Theorem \ref{th1}. 

   Suppose on the contrary that $F$ is a winning strategy for
   TWO. Now consider the game where player ONE starts the
   game by making an arbitrary move, and then also uses TWO's strategy
   $F$, while TWO uses TWO's strategy to respond to ONE.

   Consider the play $(m_1,n_1,m_2,n_2,\dots,m_k,n_k,\dots)$
   where $m_i$ is ONE's $i$-th move, $n_i$ is TWO's $i$-th move and
   $n_i=F(m_1,\dots,m_i)$ and $m_{i+1}=F(n_1,\dots,n_i)$, and
   $m_{i+1}>n_i>m_i$ for all $i$. 

   Since $F$ is a winning strategy, we have
\begin{enumerate}
\item{$\{m_j:j\in\naturals\}\in{\cal F}$ and $\{n_j:j\in\naturals\}\in
{\cal F}$ and}
\item{there is an $\ell$ such that $m_j<n_j$ for all $j>\ell$.}
\end{enumerate}

   But then, $\{m_j:j\in\naturals\}\cap \{n_j:j\in\naturals\}$ is finite,
   contradicting the fact that both sets are from the non--principal
   filter ${\cal F}$.
\end{pf}

\begin{definition} ${\cal F}$ is a rare filter if there is for
   each partition $\{I_n:n\in\naturals\}$ of $\naturals$ into disjoint
   finite sets, an $X\in{\cal F}$ such that $|X\cap I_n|\leq 1$ for each $n$.
\end{definition}

\begin{theorem}\label{rareandONE}
  The following statements are equivalent:
\begin{enumerate}
\item{ONE has a winning strategy in $\G$.}
\item{$\G$ is determined.}
\item{${\cal F}$ is not rare.}
\end{enumerate}
\end{theorem}

\begin{pf} The implication $(1)\Rightarrow(2)$ is trivial.
  
   For the implication that $(2)$ implies $(3)$,
   assume that ${\cal F}$ is rare. Since TWO does not have a
   winning strategy in this game, we consider strategies for ONE only.
   We show that ONE does not 
   have a winning strategy.  

   Consider a strategy $F$ for ONE. We may
   assume that for all $x_1<\dots<x_n$ the strategy $F$ satisfies
\[\max \{x_1,\dots,x_n\}+1<F(x_1,\dots,x_n) \mbox{ and}\]
\[F(x_1)<F(x_1,x_2)<\dots<F(x_1,\dots,x_n).\]

   Define $g$ by $g(1)=F(\emptyset)$, and
   $g(n+1)=\max\{F(j_1,\dots,j_i):j_1<\dots<j_i\leq n+1\}+g(n)$ for 
   each $n\in\naturals$. Observe that if $m<n$, then $g(m)<g(n)$.

   Put $h(n+1)=g(h(n))$ for each $n$,
   and $h(1)=F(\emptyset)+1$. Then $h$ is strictly increasing. Consider the
   partition $\{I_n:n\in\naturals\}$ where $I_n=[h(n),h(n+1))$ for $n>1$,
   and $I_1=[1,h(1))$. 

   Since ${\cal F}$ is rare, choose an $X\in{\cal F}$ such that
   $|X\cap I_n|\leq 1$ for each $n$. Enumerate $X$ in increasing order as
   $\{x_n:n\in\naturals\}$. Then choose an infinite subset $Y$ of $X$
   such that $X\setminus Y\in{\cal F}$ (for example, let $Y$ be the
   complement of a selector from ${\cal F}$ of the partition $K_1,
   K_2,\dots$ where $K_1=[1,x_1)$ and
   $K_{n+1}=[x_{n^2},x_{(n+1)^2})$). Enumerate $Y$ increasingly as 
   $\{y_1,\dots, y_n,\dots\}$. We may assume that $1<y_1$. Put
   $J_1=[1,y_1)$ and for all $n$ put $J_{n+1}=[y_n,y_{n+1})$. Since
   ${\cal F}$ is rare, we find a $Z\in{\cal F}$ such that $|Z\cap
   J_n|\leq 1$ for each $n$. Put $T=X \cap Z$. Then $T$ is also a
   selector of the family $J_1,J_2,\dots, J_n,\dots$, is in ${\cal
     F}$, and contains no endpoint of any of the $J_n$'s.

   We now have three sequences: $(n_1,n_2,\dots,n_k,\dots)$,
   $(x_{m_1},x_{m_2},\dots, x_{m_k},\dots)$, and
   $(y_{s_1},y_{s_2},\dots,y_{s_k},\dots)$, such that:
   \begin{enumerate}
     \item{$x_{m_i}\in I_{n_i}\cap J_{s_i}$ for each $i$, and}
     \item{$\max(I_{n_i})<\max(J_{s_i})<\min(I_{n_{i+1}})$ for each $i$.}
   \end{enumerate}
   But then we have for each $i$ the inequalities
\[h(n_i)\leq x_{m_i}<h(n_{i}+1)<y_{s_i}<h(n_{i+1}).\]
   
   We claim that TWO wins the play against $F$ where TWO plays
   $x_{m_1},x_{m_2},x_{m_3},\dots$. To see this, first observe that
   $F(\emptyset)<h(1)\leq h(n_1)\leq x_{m_1}<h(n_1+1)<y_{s_1}<h(n_2)$;
   thus, $F(x_{m_1})<g(x_{m_1})<g(h(n_1+1))<h(n_1+2)\leq h(n_2)\leq
   x_{n_2}$. Again applying the inequalities above and the definition
   of $h$, we see that $F(x_{m_1},x_{m_2})<x_{m_3}$, and so on.

   Thus the play $(F(\emptyset), x_{m_1}, F(x_{m_1}), x_{m_2},
   F(x_{m_1},x_{m_2}), x_{m_3},\dots)$ has the properties that the set
   of moves by TWO is in ${\cal F}$, and as a sequence eventually
   dominates the sequence of moves by ONE. It follows that $F$
   is not a winning strategy of ONE. This completes the proof of
   $\neg(3)\Rightarrow\neg(2)$.

   Next we show that $(3)$ implies $(1)$: Let ${\cal F}$ be a 
   non--rare filter, and choose a 
   partition $\{I_n:n\in \naturals\}$ of $\naturals$ into disjoint finite sets
   such that each element of ${\cal F}$ meets infinitely many of the
   $I_n$-s in more than one point. Player ONE's strategy $F$ will be a simple
   $1$--tactic (i.e., it depends only on the most recent move of the
   opponent). Define 
   $F$ as follows: Let $k$ be given. Fix $n$ so that $k\in I_n$ and put
   $F(k)=(\max\cup_{j\leq n+k}I_j)+ n + k + 1$.

   To see that $F$ is winning for ONE, suppose that
   $(m_1,n_1,m_2,n_2,\dots,m_k,n_k,\dots)$ is an $F$--play. Then
   $m_{k+1}=F(n_k)$ for each $k$. For each $k$ let $j_k$ be such that
   $n_k\in I_{j_k}$. By the definition of $F$ we see that if for all but
   finitely many $k$ we have $m_k<n_k$, then for all but finitely many
   $k$ we have $m_{k+1}<j_{k+1}<m_{k+2}<j_{k+2}$; thus
   $X=\{n_k:k\in\naturals\}$ meets all but finitely many $I_n$ in at most
   one point. But then $X\not\in{\cal F}$.

   Thus, either TWO's sequence of moves enumerates a set in the filter
   but does not eventually dominate ONE's sequence of moves, or else
   TWO's sequence of moves eventually dominates ONE's sequence of
   moves, but does not enumerate an element of ${\cal F}$. In either
   event ONE wins.
\end{pf}

\section{Remarks.}

   Both of our games can be coded as Gale--Stewart games in such a way
   that if the filter is projective then the corresponding
   Gale--Stewart game is projective. Assume the Axiom of Dependent
   Choices. Determinacy hypotheses such
   as the Axiom of Determinacy imply that our games are
   determined, and thus that all filters are meager and
   non--rare. (Projective Determinacy would imply that our games
   are determined for all projective filters, and thus that all
   projective filters are meager (hence non--rare).) But a much
   weaker hypothesis, namely that every set of reals has the property
   of Baire, already implies the determinacy of our games. In the case
   of rare filters this follows because Talagrand's theorem
   implies that these do not have the   property of Baire. 

   In $ZFC$, determinacy of the games $\G$ is weaker than determinacy of $\HH$.
   This can be seen as follows: In Theorem 5.1 of his paper \cite{K}, Kunen
   shows that in the random real model there are no rare
   ultrafilters -- and thus no rare filters. In this model $\G$
   is determined. But there are always non--Baire filters, since every
   non--principal ultrafilter is like that.

   It is well--known that in the presence of the Continuum Hypothesis
   or Martin's Axiom, there are rare filters and thus
   undetermined instances of our games. It can also be insured that
   these undetermined filters are (are not), $P$-filters, and so on.

\end{document}